\documentclass[a4paper,10pt]{article}
\usepackage{amsmath,amssymb,amsthm,latexsym,graphicx}
\pdfoutput=1

\newtheorem{lemma}{Lemma}

\newtheorem{theorem}{Theorem}
\newtheorem{prop}{Proposition}

\numberwithin{equation}{section}

\title{\sc Transferring Fourier Multipliers from $\mathbb{R}^n$ to Compact Lie Groups}
\author{{David G. Maher}\footnote{
$\langle$DavidGMaher@yahoo.com.au$\rangle$ }}

\begin{document}

\maketitle

\newcommand{\g}{\mathfrak{g}}
\newcommand{\kg}{\mathfrak{k}}
\newcommand{\tg}{\mathfrak{t}}
\newcommand{\ug}{\mathfrak{u}}
\newcommand{\p}{\mathfrak{p}}
\newcommand{\X}{\mathfrak{X}}
\newcommand{\af}{\mathfrak{a}}
\newcommand{\R}{\mathbb{R}}
\newcommand{\C}{\mathbb{C}}
\newcommand{\Z}{\mathbb{Z}}
\newcommand{\N}{\mathbb{N}}
\newcommand{\E}{\mathbb{E}}
\newcommand{\bbP}{\mathbb{P}}
\newcommand{\T}{\mathbb{T}}
\newcommand{\F}{\mathcal{F}}
\newcommand{\A}{\mathcal{A}}
\newcommand{\x}{\mathbf{x}}
\newcommand{\y}{\mathbf{y}}
\newcommand{\ad}{\mathrm{ad}}
\newcommand{\Ad}{\mathrm{Ad}}
\newcommand{\Exp}{\mathrm{Exp}}
\newcommand{\grad}{\mathrm{grad} \,}
\newcommand{\dd}{\mathrm{\, d}}

\begin{abstract}
This paper uses the wrapping map of Dooley and Wildberger to prove $L^p$ boundedness of multipliers on compact Lie groups by transferring the estimate from $\mathbb{R}^n$. This improves the bounds in several cases, and simplifies the proofs of others.\\
\end{abstract}

\underline{\bf Contents}\\
\phantom{abcde} 1. Introduction\\
\phantom{abcde} 2. Notation and Preliminaries\\
\phantom{abcde} 3. The Wrapping Map and Multipliers\\
\phantom{abcde} 4. Main Results: $L^p-L^q$ bounds\\
\phantom{abcde} 5. An Application: Convergence of Fourier Series\\
\phantom{abcde} 6. Further directions\\
\phantom{abcde} 7. References\\

\section{Introduction}

Let $G$ be a compact connected Lie group and $\g$ its associate Lie algebra.

There are many papers dealing with the $L^p$ boundedness of
multipliers on compact Lie groups, including \cite{ST}, \cite{STR},
\cite{W}.  These papers have made use of one or both of the
following techniques: the first is to try and transfer the analysis
from the curved manifold to $\R^n$.  The main drawback here is that
this can only occur within a neighbourhood of where the exponential
map is injective.\\

The second technique is to transfer the analysis directly to a
maximal torus $T \subseteq G$ to appeal to existing (abelian)
Fourier analysis.  This has involved studying the integrability of
negative powers of the Weyl function to transfer estimates on $G$ to
estimates on $T$.  Since this function vanishes at the origin, this creates many technical problems.\\

The wrapping map $\Phi$ introduced in \cite{DW1} is a homomorphism of central measures
or distributions under (abelian) convolution on $\g$ (viewed as $\R^n$)
to their corresponding counterparts on $G$, utilising the whole of the lift of the exponential map,
not just its restriction to a fundamental domain.  That is,
$$
\Phi(\mu *_{\g} \nu) = \Phi(\mu) *_G \Phi(\nu)
$$

As any Fourier multiplier can be written as a convolution with the associated kernel operator, the wrapping map provides a natural mechanism to transfer results concerning Fourier multipliers from $\g \cong \R^n$ to $G$.\\


\section{Notation and Preliminaries}

Let $G$ be a compact connected Lie group, $\g$ its Lie algebra, $T$ a maximal torus and $\tg$ the Lie algebra of
$T$.  Let $n$ be the dimension of $G$, and $l$ the dimension of $T$ (also known as the rank of $G$).  Let $B( \cdot , \cdot )$ be the Killing form on $\g$, with $\g^*$ and $\tg^*$ the respective duals of $\g$ and $\tg$ with respect to the Killing form.\\

We denote by $\Sigma^+$ the set of positive roots $\{ \alpha_1, \dots, \alpha_k \}$, with $k = \bigl|
\Sigma^+ \bigr|$, and thus we have $n = l + 2k$.  Let $W$ denote the Weyl group, $\tg^*_+$ the positive Weyl chamber, and let $\rho = \tfrac 12 \sum_{\alpha \in \Sigma^+} \alpha$.  Let $\Lambda \subseteq \tg^*$ denote the set of integral weights, and $\Lambda^+ = \Lambda \cap \tg^*_+$  the set of positive integral weights.\\

Every irreducible representation $\pi \in \widehat{G}$ is associated
to a unique highest weight $\lambda \in \lambda^+$.  If $\chi_\pi =
\chi_\lambda$ is the character of this representation, the
Kirillov's character formula is given by
$$
j(X) \chi_\lambda (\exp X) = \int_{\mathcal{O}_{\lambda + \rho}}
e^{i\beta(X)} d\mu_{\lambda + \rho} (\beta), \; \; \; \text{for all}
\; \; \; X \in \g
$$
where $\mathcal{O}_{\lambda + \rho}$ is the co-adjoint orbit through
$\lambda + \rho \in \tg^*_+$, and $\mu_{\lambda + \rho}$ is the
Liouville measure on $\mathcal{O}_{\lambda + \rho}$ with total mass
$d_\lambda = d_\pi = \dim \pi$, and where $j(X)$ is is the analytic square root of the Jacobian of
$\exp$ with $j(0) = 1$, given by
$$
j(X) = \prod_{\alpha \in \Sigma^+} \frac{\sin \alpha (X)/2}{\alpha
(X)/2}
$$
That is, for $f \in C_c^\infty (G)$:
$$
\int_{\g} f(\exp X) |j(X)|^2 dX = \int_{G} f(g) dg:
$$
where we have normalised the Haar measure $dg$ on $G$ to have total
mass 1, with subsequent normalisation of the Lebesgue measure $dX$ on $\g$. Similarly, we normalise the Haar measure and $dt$ on $T$ and the Lebesgue measure $dH$ on $\tg$ so that we have for $f \in C_c^\infty (T)$:
$$
\int_T f(\exp H) dH = \int_{\tg} f(t) dt
$$

Integrals over $G$ of central functions may be computed by integrals
on $T$ by the Weyl integral formula:
$$
\int_G f(g)dg = \frac{1}{|W|} \int_T  f(t) |\triangle (t)|^2 dt, \; \text{ where } \; \triangle (\exp H) = \prod_{\alpha \in \Sigma_+} 2 \sin \alpha(H)/2
$$

Correspondingly, we also have a similar integral formula for an Ad-invariant function of compact support $\varphi$ on $\g$:
$$
\int_\g \varphi(X)dX = \int_{t_+} \prod_{\alpha \in \Sigma_+}
\alpha(H)^2 \varphi(H) dH.
$$

We also denote by $L^p (\g)$, $L^p (\tg)$, $L^p (G)$ and $L^p (T)$
the set of $p$-integrable functions on $\g$, $\tg$, $G$ and $T$,
respectively, and $\| \cdot \|_p$ the $p$-norm when the space
concerned is unambiguous.  We also denote by $L^p_I (\g)$ and $L^p_I
(G)$ the respective sets of $G$-invariant $L^p$ functions.\\

Let $\mu^\wedge$ denote the Fourier transform of $\mu$, with the convention
$$
\mu^\wedge (\xi) = \int_{\R^n} \mu (\x) e^{i \x \cdot \xi} d\x .
$$
when $\mu \in L^1 (\R^n)$.  Correspondingly, Let $\pi \in \widehat{G}$ have highest weight $\lambda \in \Lambda^+$.  The Fourier transform of a central function, $\nu$, is a multiple $c_\pi I_\pi$ of the identity, where
$$
c_\pi = \frac{1}{d_\pi} \int_G \nu (g) \chi_\lambda (g) dg
$$
when $\mu \in L^1_I (G)$, such that we may write the Fourier series
$$
\mu (g) = \sum_{\lambda \in \Lambda^+} c_\lambda d_\lambda \chi_\lambda (g) .
$$

We also employ the notation $X \lesssim Y$ to refer to the situation $X \leq CY$, where $C$ is an unspecified constant.\\

\underline{\bf The Wrapping Map}\\

Assume that $\nu$ is a distribution of compact support on $\g$, or
$j\nu \in L^1(\g)$.  We define the {\bf wrapping map}, $\Phi$ by
\begin{equation}
\langle \Phi(\nu), f \rangle = \langle \nu, j \tilde{f} \rangle
\end{equation}
where $f \in C^\infty (G)$, $\tilde{f} = f \circ \exp$.  We
call $\Phi(\nu)$ the {\bf wrap} of $\nu$.  The principal result is
the {\bf wrapping formula} (\cite{DW2}, Thm. 2) given by
\begin{equation}\label{wrapform}
\Phi(\mu *_{\g} \nu) = \Phi(\mu) *_G \Phi(\nu)
\end{equation}

What (\ref{wrapform}) shows us is that problems of convolution of
central measures or distributions on a (non-abelian) compact Lie
group can be ``transferred" to Euclidean convolution of Ad-invariant
distributions on $\g$.  Thus, since a Fourier multiplier operator
may be regarded as a convolution operator, a natural question arises
as to how we may ``wrap" a Fourier multiplier operator.  To explicitly show this connection, we must consider the formulation for $\Phi$.

\begin{prop} We have the following formulations for $\Phi$ for $G$-invariant distributions of compact support, firstly given as summations over the integer lattice:

\begin{equation}\label{PhiWrap}
\Phi(\mu) (\exp \, H) = \sum_{\gamma \in \Gamma}
\frac{\mu}{j} (H + \gamma) \; \; \; \forall H \in \tg.
\end{equation}

\begin{equation}\label{PhiWrapj}
\Phi(j \mu) (\exp \, H) = \sum_{\gamma \in \Gamma} \mu (H +
\gamma) \; \; \; \forall H \in \tg.
\end{equation}

and as the Fourier series:

\begin{equation}\label{PhiFS}
\Phi (\mu)(\exp H) = \sum_{\lambda \in \Lambda^+} d_\lambda
\mu^\wedge (\lambda + \rho) \chi_\lambda (\exp H)
\end{equation}

\begin{equation}\label{PhiFSj}
\Phi (j \mu)(\exp H) = \sum_{\lambda \in \Lambda^+} d_\lambda
\int_G \mu^\wedge (\lambda + \rho -
\ad (g)(\rho)) \dd g \, \chi_\lambda (\exp H)
\end{equation}

where the values are given on $H \in \tg$

\end{prop}

(\ref{PhiWrapj}) is calculated in \cite{DW2}, and (\ref{PhiWrap}) follows naturally.  In \cite{DW2} it is shown that the Fourier transform of $\Phi(\mu)$ at $\pi$ is given by
$$
c_\pi = (\Phi(\mu)) ^\wedge (\lambda + \rho) = \mu^\wedge (\lambda + \rho)
$$
by applying Kirillov's character formula, that is:
$$
c_\pi = \frac{1}{d_\pi} \langle \Phi(\mu) , \chi_\lambda \rangle =
\frac{1}{d_\pi} \langle \mu , j \tilde{\chi}_\lambda \rangle = \langle \mu , \int_G e^{ig(\lambda + \rho)(\cdot)} dg
\rangle
$$

By the $G$-invariance of $\mu$ this is
\begin{equation}\label{cpi}
c_\pi = \langle \mu , e^{i(\lambda + \rho)(\cdot)} \rangle =
\mu^\wedge (\lambda + \rho)
\end{equation}

which yields (\ref{PhiFS}).  Note that (\ref{wrapform}) also follows from (\ref{cpi}), since
\begin{align*}
(\Phi(\mu * \nu))^\wedge (\pi) & = (\mu * \nu)^\wedge (\lambda + \rho) I_\pi\\
& = \mu ^\wedge (\lambda + \rho) \nu ^\wedge (\lambda + \rho) I_\pi\\
& = (\Phi(\mu)) ^\wedge (\lambda + \rho) (\Phi(\nu)) ^\wedge (\lambda + \rho) I_\pi\\
& = (\Phi(\mu) * \Phi(\nu))^\wedge (\pi)
\end{align*}

We now have an analogue of (\ref{cpi}) for Ad-invariant
distributions with compact support of the form $j \mu$:
\begin{lemma}\label{PhiHatj}  Let $\mu$ be an Ad-invariant distribution of compact support on $\g$.  Then the Fourier transform of $\Phi(\mu)$ at $\pi$ is a multiple $c_\pi I_\pi$ of the identity, where
$$
c_\pi = (\Phi(\mu)) ^\wedge (\lambda + \rho) = \int_G \mu^\wedge (\lambda + \rho -
\ad (g)(\rho)) \dd g
$$
\end{lemma}

{\bf Proof:} From (\ref{cpi}) we have that as a multiple of
the identity:
$$
\Phi(\mu)^\wedge (\pi_\lambda) = \mu^\wedge (\lambda + \rho)
$$

Thus,
\begin{align*}
\Phi(j \mu)^\wedge (\pi_\lambda) & = (j \mu)^\wedge
(\lambda + \rho)\\
& = (j^\wedge * \mu^\wedge) (\lambda + \rho)\\
& = \int_G \mu^\wedge (\lambda + \rho - \ad (g)(\rho)) \dd g
\phantom{abcde} \square
\end{align*}

Reconciling (\ref{PhiWrap}) with (\ref{PhiFS}), or (\ref{PhiWrapj}) with (\ref{PhiFSj}), gives the Poisson summation
formula for compact Lie groups.\\

In the next section we will consider these two function spaces in
the treatment of multipliers.

\section{The Wrapping Map and Multipliers}

Consider the Fourier multiplier operator, $T$, on $\L^p$ by
$$
\widehat {T f} (\xi) = m(\xi) \hat{f} (\xi)
$$

where $m(\xi)$ is the {\bf symbol} of the operator $T$.  By the standard kernel theorems, every Fourier
multiplier $T$ maybe expressed as a convolution operator with kernel $K$, that is: $T f = K * f$.\\

On $\g \cong \R^n$, a multiplier $T_\psi$ takes the form:
$$
(T_\psi f)(x) = \frac{1}{(2\pi)^n} \int_{\g} \psi(\xi) f(\xi) e^{-i
B(\xi, x)} d\xi
$$
with kernel
$$
K_\psi (x) = \int_{\g} \psi(\xi) e^{i B(\xi, x)} d\xi
$$
that is
$$
\widehat{K_\psi} (\xi) = \psi(\xi)
$$

On $G$, a multiplier $T_\Psi$ takes the form:
$$
(T_\Psi f)(g) = \sum_{\lambda \in \Lambda^+} \Psi(\lambda) f(\lambda)
\chi_\lambda (g)
$$
with kernel
$$
K_\Psi (g) = \sum_{\lambda \in \Lambda^+} \Psi(\lambda) \chi_\lambda (g)
$$
that is
$$
\widehat{K_\Psi} (\lambda) = \Psi(\lambda)
$$

Wrapping the multiplier $T_\psi$ from $\g \cong \R^n$ to $G$ we have
$$
\Phi(T_\psi f) = \Phi(K_\psi * f) = \Phi(K_\psi) * \Phi(f)
$$

To compute the kernel $\Phi(K_\psi)$, we see that its Fourier
coefficients are
$$
\Phi^\wedge (K_\psi)(\pi_\lambda) = \Psi (\lambda) = K_\psi^\wedge
(\lambda)
$$

Thus, we need to compute $\Psi (\lambda) = K_\psi^\wedge (\lambda)$.  In
light of (\ref{cpi}), we see that
\begin{equation}\label{m1}
\Psi (\lambda) = \psi (\lambda + \rho)
\end{equation}

Alternatively, if we wish to consider $\Phi(j \cdot K_\psi)$, we see by
Lemma \ref{PhiHatj} that
\begin{equation}\label{m2}
\Psi (\lambda) = \int_G \psi (\lambda + \rho - \ad (g)(\rho)) \dd g
\end{equation}

These two expressions \ref{m1} and \ref{m2} for the multiplier $\Psi$
on $G$ constitute to two forms considered in \cite{STR}, where they
are referred to as ``(*)'' and ``(**)'', respectively.\\

Thus, to obtain bounds for multipliers on $L^p_I(G)$, we now need to
consider bound of the form $\| \Phi (\nu) \|_q \lesssim \| \nu \|_p$
and $\| \Phi (j \cdot \nu) \|_q \lesssim \| \nu \|_p$

\section{Main Results: $L^p-L^q$ bounds}

In this section we consider the following problem: Suppose $\nu \in
L^p(\g)$.  For what $q$ is the wrap of $\nu$ in $L^q (G)$?  That is,
for what $p$ and $q$ do we have
$$
\| \Phi (\nu) \|_q \lesssim \| \nu \|_p
$$

Furthermore, for what $p$ and $q$ do we have
$$
\| \Phi (j \nu) \|_q \lesssim \| \nu \|_p
$$

These will then be applied to the $L^p-L^q$ bounds of multipliers.\\

We firstly have the following key result:
\begin{lemma}\label{jLp} $j \in \L^p (\g)$ for $p > 2n / (n - l)$.
\end{lemma}

{\bf Proof:}  Firstly, note that $j$ is bounded at $0 \in \g$.  Denote a neighbourhood $B$ of $0 \in \g$, with $\int_{B} |j(X) |^p \dd X = C$ being a finite quantity.
\begin{align*}
\| j(X) \|^p_p & = \int_{\g} | j(X) |^p \dd X \\
& = \int_{\g \setminus B} |j(X) |^p \dd X + C\\
& \leq \int_{\g \setminus B} \; \Bigl| \prod_{\alpha \in \Sigma_+} \alpha (X/2) \Bigr|^{-p} \dd X + C\\
& = \int_{\tg_+ \setminus B} \; \Bigl| \prod_{\alpha \in \Sigma_+} \alpha (H/2) \Bigr|^{2-p} \dd H + C\\
\end{align*}

The function $\prod_{\alpha \in \Sigma_+} \alpha (H/2)$ is a
homogeneous polynomial of degree $k$ on a space of dimension $l$,
having positive co-efficients, and thus will be bounded provided $k
(2-p) < l$.  That is, $p > 2n / (n - l)$. $\phantom{abcde}
\square$\\


As a consequence of this bound on $j$, we have the following bounds on $p$ and $q$ for which $\| j f \|_q \lesssim \| f \|_p$:
\begin{lemma}\label{jf1}  We have

\phantom{abc} a) $\| j f \|_p \leq \| f \|_p$ for $1 < p < \infty$.

\phantom{abc} b) $\| j f \|_p \leq \| f \|_\infty$ for $p > 2n / (n - l)$.

\phantom{abc} c) $\| j f \|_1 \leq \| f \|_q$ for $q < 2n / (n + l)$\\
\end{lemma}

{\bf Proof:} These follow from the Hausdorf-Young, Young, and H\"{o}lder inequalities, and are thus best possible.\\

{\bf Proof of a):}  Let $p$ and $q$ be conjugate exponents.  We have
\begin{align*}
\| j f \|_p & \leq \| \widehat{j f} \|_q \tag{by Hausdorf-Young}\\
& = \| \hat{j} * \hat{f} \|_q\\
& \leq \| \hat{j} \|_1 \| \hat{f} \|_q \tag{by Young}\\
& = \| \mu_\rho \|_1 \| f \|_p\\
& = \| f \|_p
\end{align*}
since the mass of the Liouville measure $\mu_\rho$ is equal to $\dim \pi_0 = 1$.

{\bf Proof of b):}  Let $p$ and $q$ be conjugate exponents.  We have
\begin{align*}
\| j f \|_p & \leq \| \widehat{j f} \|_q \tag{by Hausdorf-Young}\\
& = \| \hat{j} * \hat{f} \|_q\\
& \leq \| \hat{j} \|_q \| \hat{f} \|_1 \tag{by Young}\\
& = \| j \|_p \| f \|_{\infty}\\
\end{align*}
which is bounded so long as $p > 2n / (n - l)$, since we have $j \in
\L^p (\g)$ by Lemma \ref{jLp}.

{\bf Proof of c):}  Let $p$ and $q$ be conjugate exponents.  By H\"{o}lder
we have:
$$
\| j f \|_1 \leq \| j \|_p \| f \|_q
$$

By Lemma \ref{jLp}, $j \in \L^p (\g)$ for $p > 2n / (n - l)$. Hence,
the result follows if $f \in \L^q (\g)$ for $q < 2n / (n + l)$. $\phantom{abcde} \square$\\

Regarding $L^p$ bounds for $\Phi$, we have the following from \cite{DW2}, which is almost obvious from the definition: If $\nu \in L^1 (\g)$, then we have:
\begin{equation}\label{nuL1}
\| \Phi (\nu) \|_1 \leq \| \nu \|_1
\end{equation}
\begin{equation}\label{jnuL1}
\| \Phi (j \nu) \|_1 \leq \| \nu \|_1
\end{equation}

We now prove a more general $L^p$ bound:

\begin{theorem}\label{jf4} We have the bounds
\begin{equation}
\| \Phi (j \mu) \|_p \lesssim \| \mu \|_p , , \; \; 1 \leq p \leq \infty
\end{equation}

\begin{equation}
\| \Phi (\mu) \|_p \lesssim \| \mu \|_p \, , \; \; \; p \geq 2
\end{equation}
\end{theorem}

{\bf Proof:} Let $\tg_\Gamma \subseteq \tg$ be a fundamental domain for $\Gamma$ in $\tg$.  We have:
\begin{align*}
\| \Phi (j \mu) \|_p^p &= \int_G | \Phi (j \mu)(g) |^p \dd g\\
&= \frac{1}{|W|} \int_T | \Phi (j \mu)(t) |^p | \triangle(t) |^2 \dd t\\
&= \frac{1}{|W|} \int_{\tg_\Gamma} | \sum_{\gamma \in \Gamma} \mu(H + \gamma) |^p | \triangle(\exp H) |^2 \dd H \\
&\leq \frac{1}{|W|} \int_{\tg} | \mu(H) |^p | \triangle(\exp H) |^2 \dd H \\
&= \int_{\tg^+} | \mu(H) |^p | j(H) |^2 \Bigl| \prod_{\alpha \in \Sigma_+} \alpha (H/2) \Bigr|^2 \dd H\\
&= \int_{\g} | \mu(X) |^p | j(X) |^2 \dd X \\
&\leq \int_{\g} | \mu(X) |^p \dd X = \| \mu(X) \|_p^p
\end{align*}

Similarly for $\Phi (\mu)$, we arrive at
\begin{align*}
\| \Phi (\mu) \|_p^p &= \int_{\g} | \tfrac{\mu}{j}(X) |^p | j(X) |^2 \dd X \\
&= \int_{\g} | \mu (X) |^p | j(X) |^{2-p} \dd X \\
&\leq \| \mu \|_p^p \| | j(X) |^{2-p} \|_{\infty}\\
&\leq \| \mu \|_p^p
\end{align*}
as long as $p \geq 2$. $\phantom{abcde} \square$\\

{\bf Remark:} Compare this to \cite{STR} Lemma 3, which asserts these bounds as long as
$$
\sum_{\gamma \in \Gamma} |j(H + 2\pi \gamma)|^{np'-2}
$$
is bounded, which is only explicitly addressed for the case of $SU(2)$.  Here, $n = 1$ is the case for $\Phi (\mu)$, and $n = 2$ is the case for $\Phi (j \mu)$.  However, there appears to be an error in this calculation.\\

We can use this to prove further $L^p-L^q$ for $\Phi (j \mu)$ only.

\begin{theorem}\label{jf3} We have the $L^q-L^1$ bound
\begin{equation}
\| \Phi (j \mu) \|_1 \lesssim \| \mu \|_q \, , \; \; \; q < 2n / (n + l)
\end{equation}
and the $L^\infty-L^q$ bound
\begin{equation}
\| \Phi (j \mu) \|_q \lesssim \| \mu \|_\infty \, , \; \; \; q > 2n / (n - l)
\end{equation}
\end{theorem}

{\bf Proof:} Let $\nu = j \mu$.  From, (\ref{nuL1}) we have that
$$
\| \Phi (\nu) \|_1 \lesssim \| \nu \|_1
$$
and so by setting $\nu = j \mu$, we have
$$
\| \Phi (j \mu) \|_1 \lesssim \| j \mu \|_1
$$

By Lemma \ref{jf1} c), it is only possible to have
$$
\| \Phi (j \mu) \|_1 \lesssim \| j \mu \|_1 \leq \| \mu \|_q
$$
for $q < 2n / (n + l)$.  Similarly, By lemma \ref{jf1} b), it is only possible to have
$$
\| \Phi (j \mu) \|_q \lesssim \| j \mu \|_q \leq \| \mu \|_\infty
$$
for $q > 2n / (n - l)$.  $\phantom{abcde} \square$\\

We concluded this section with an additional proof of the $L^2$ bound for $\Phi (j \nu)$ only.  This uses the fact that the Fourier transform of $j$ leads to a region around each point $\lambda \in\Lambda^+$ to tessellate over $\tg^*$.  We then apply Parseval.
\begin{theorem}
$$
\bigl\| \Phi (j \nu) \bigr\|_2 \leq \bigl\| \nu \bigl\|_2
$$
\end{theorem}

{\bf Proof:} From Lemma \ref{PhiHatj} we have that
$$
d_\lambda \Phi(j \nu)^\wedge (\pi_\lambda) = \int_G \nu^\wedge
(\lambda + \rho - \ad (g)(\rho)) \dd g
$$

Let $Q$ be the convex hull of $\{ w \rho \; | \; w \in W \}$ and $E$
a bounded function on $Q$.  We have,
\begin{align*}
\bigl\| \Phi (j \nu) \bigr\|^2_2 & = \sum_{\lambda \in
\Lambda^+} \biggl| \int_G \nu^\wedge (\lambda + \rho - \ad
(g)(\rho)) \dd g \biggr|^2\\
& \leq \sum_{\lambda \in \Lambda^+} \biggl| \int_Q \biggl(
\prod_{\alpha \in \Sigma^+} \frac{\partial}{\partial \alpha} \biggr)
\nu^\wedge (\lambda + \rho - \xi) E(\xi) \dd \xi \biggr|^2\\
& \lesssim \int_{\tg^*} \, \biggl| \biggl(\prod_{\alpha \in \Sigma^+} \frac{\partial}{\partial \alpha} \biggr) \nu (\xi) \biggr|^2 \dd \xi\\
& = \int_{\tg} \, \biggl| \prod_{\alpha \in \Sigma^+} \alpha (H) \nu (H) \biggr|^2 \dd H\\
& \lesssim \int_{\g} \, | \nu (X) |^2 \dd X\\
& = \bigl\| \nu \bigl\|_2^2 \phantom{abcde} \square
\end{align*}

%
%
%

\section{An Application: Convergence of Fourier series}

We firstly require some notation, which we adopt from \cite{ST}. Let
$R$ be a Weyl-invariant, closed, convex polyhedral subset of $\tg$
which contains the origin.  Let $t R = \{ tX | X \in R\}$.  As $t$
ranges over $[1,\infty)$, $(tR) \cap \tg^*_+$ generates only a
countable number of distinct sets, which is denoted by $\{ R_N \}$.  Define the Weyl-invariant polygonal operator on $\g \cong \R^n$ by:
$$
K_t \mu (x) = \int_{tR} \exp(-isx) (\hat\mu (s)) ds
$$

We define the $N$th partial sum of the Fourier series of a function
$ f \in L^p(G), \; 1 \leq p \leq \infty$, as

$$
S_N f(g) = \sum_{\lambda \in R_N} d_\lambda \chi_\lambda * f(g)
$$

Then we have

\begin{theorem}  Suppose $R$ is a regular polyhedron, then if $f \in L^p_I
(G), \; p > 2n / (n + l)$, then $S_N f$ converges to $f$ almost everywhere.
\end{theorem}

{\bf Proof:} Since $\Phi (K_t) = S_N$, and from \cite{F} the Fourier transform converges for any polygonal region containing the origin for $1 < p < \infty$, and therefore from section 4 $\Phi (K_t)$ is bounded for $1 \leq p \leq \infty. \phantom{abcde} \square$

\section{Further Directions}

The wrapping map devised in \cite{DW1} has already been employed as a transference method in analysis.  These could be further extended using the estimates (or the idea of transferring from the tangent space to a curved space) in this paper.  For example:
\begin{itemize}
\item In \cite{DW2} the concept of a `modulator' was introduced.  In particular, this allows one to consider the class of modulators that `wraps' to a (generalised) character.  Thus, one approach to considering the norms of characters on $G$ would be to consider the norm of the modulator on $\g$.
\item In \cite{M1} the wrapping map was used to `wrap' not only heat kernels but also Brownian motion from $\g \cong \R^n$ to $G$.  The estimates in this paper could be used to compute $L^p$ norms for solutions for these and other partial differential equations on compact Lie groups, as well as Brownian motions or other stochastic processes.
\item In \cite{M2} the wrapping map was extended to complex Lie groups and compact symmetric spaces.  $L^p$ bounds on these spaces could again be computed by considering the $L^p$ bounds on their respective Lie algebra and tangent space.
\end{itemize}

We will consider these in future work.

\bigskip

%
%

Email: {\tt DavidGMaher@yahoo.com.au}

\end{document}